\documentclass{amsart}
\usepackage{verbatim,amsfonts,color,mathrsfs,stmaryrd,graphicx, mathdots}
\usepackage{xcolor}
\usepackage{tikz}
\usepackage{tikz-cd}
\usepackage{tkz-euclide}
\usetikzlibrary{patterns}   
\usetikzlibrary{positioning}
\usetikzlibrary{decorations.pathreplacing,calligraphy}

\usepackage{pgfplots}

\pgfplotsset{compat=1.6}
\pgfplotsset{soldot/.style={color=blue,only marks,mark=*}} 
\pgfplotsset{holdot/.style={color=blue,fill=white,only marks,mark=*}}

\usepackage{diagbox} 
\usepackage{tikz, pgf, pgfplots}
\pgfplotsset{compat=newest}

\usetikzlibrary[patterns]
\usetikzlibrary{arrows, arrows.meta}
\usetikzlibrary{calc} 
\def\centerarc[#1](#2)(#3:#4:#5)%
    {\draw[#1] ($(#2)+({#5*cos(#3)},{#5*sin(#3)})$) arc (#3:#4:#5);}

\usetikzlibrary{intersections,through,backgrounds} 
\usetikzlibrary{shapes.geometric} 

\usepackage{enumitem}
\usepackage{amssymb}

\theoremstyle{plain}

\theoremstyle{definition}

\theoremstyle{remark}

\usepackage{hyperref}

\begin{document}

 
\title[Overview of C.~Series' {\em The modular surface and continued fractions}]{An elegant model of the geodesic flow on the modular surface}
\author{Pierre Arnoux}
\address{Universit\'e d'Aix-Marseille\\ 
Institut de Math\'ematiques de Marseille\\
I2M - UMR 7373\\
13453 Marseille, France}
\email{arnoux@iml.univ-mrs.fr}

\author{Thomas A. Schmidt}
\address{Department of Mathematics\\ 
              Oregon State University\\
              Corvallis, OR 97331\\
              USA}
\email{toms@math.orst.edu}
\date{3 May 2026}

\thanks{It is a pleasure to thank Pascal Hubert for discussion,  and the referees for suggestions and corrections.}
 

\begin{abstract}  Caroline Series' [{\em The modular surface and continued fractions}, J. Lond. Math. Soc. (2), {\bf 31}, no.~1, (1985), 69--80] gives a clear framework linking, in a deceptively simple way, the dynamics of the geodesic flow on the modular surface  with the dynamics of the regular continued fraction, through a well-chosen symbolic coding.     It has been called {\em required reading}  for those interested in the symbolic dynamics of geodesic flows,  and has had consequences in   symbolic dynamics, ergodic theory, hyperbolic geometry, and continued fraction theory.    In this overview, we give an indication of why this is so,  sketch some of the history related to the paper, and also point to some later works. 
\end{abstract}

\maketitle

\tableofcontents

\section{Introduction} 


The 1985 article \cite{Series} of Caroline Series  
  gives  
a clear framework linking, in a deceptively simple way, the dynamics of a   three dimensional  geometric dynamical system with the dynamics of a  one dimensional system. To wit, the paper  connects the
 geodesic flow on the modular surface to the dynamics of the usual continued fraction, through a  symbolic coding    that is  elegantly simple.  The paper itself is  short, precise, and accessible. \footnote{ On a  
 personal note, one of us (TS) learned of \cite{Series} in what now is perhaps viewed as a quaint manner.   He remains very grateful to the person who left a reprint of the paper in a common room of the department where he was a graduate student.}

The connection of this geodesic flow to continued fractions lies at the nexus of dynamics, geometry and number theory and hence ``{\em The modular surface and continued fractions}" has had an impact well beyond the studies of merely the modular surface or of continued fractions.  The very elegance of the Farey tesselation model and the coding given in the paper have naturally led  many researchers  to extend this model to other situations. 

  We end this introduction by giving the sense of the main result of \cite{Series},  preceded by  a terse review of background.   In the following two sections,  we give more detail:  a sketch of the history leading up to \cite{Series}  and a  brief  summary of the paper itself.   In the ensuing sections, we mention alternate methods for coding the geodesics of the modular surface, and then point to a selection of  various and diverse works which cite \cite{Series}.   The concluding paragraph includes a list of some introductory texts on materials related to \cite{Series}.

Since at least Poincar\'e's  insights on the many bodies problem,  mathematicians have sought to describe the topology of dynamical systems' orbits.   Naturally occurring dynamical systems  --- such as the collection of orbits of the bodies in our solar system ---  consist of objects following  energy minimizing paths.   The mathematical notion of `geodesic',  locally length minimizing paths, abstracts this idea.    In 1898, Hadamard \cite{Hadamard} gave symbolic codings of geodesics, doing so in the setting of noncompact surfaces of negative curvature with infinite funnels.  Very roughly,  he associated to a geodesic avoiding the funnels an infinite cutting sequence given by its encounters with a marked set of geodesic segments on the surface.   This allowed him to prove the density of the periodic geodesics in the set of all  bounded geodesics.    From many a viewpoint, the modular surface is the quintessential noncompact surface of negative curvature, but it has two cone points and a cusp instead of funnels and hence Hadamard did not consider it.

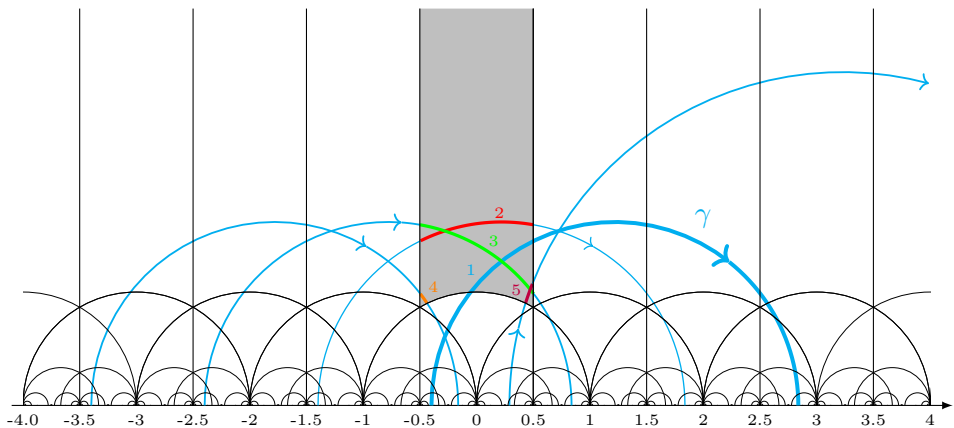
\begin{figure}[h!]
%
\pgfmathsetmacro{\myxlow}{-4}
\pgfmathsetmacro{\myxhigh}{4}
\pgfmathsetmacro{\myiterations}{12}

\begin{tikzpicture}[scale=0.5,x=3.0cm,y=3.0cm]
    \draw[-latex](\myxlow-0.1,0) -- (\myxhigh+0.2,0);
    \pgfmathsetmacro{\succofmyxlow}{\myxlow+0.5}
   \foreach \x in {\myxlow,\succofmyxlow,...,\myxhigh}
   \draw (\x,0) -- (\x,0) node[below,font=\tiny] {\x};
\node[cyan] at (2, 1.65)  {\large $\gamma$}; 
 \node[cyan] at (-0.05, 1.2)  {\tiny $1$}; 
\draw [shift={(1.22,0)},line width=1.5pt,color=cyan, ->]  plot[domain=pi:pi/3.5,variable=\t]({1.618*cos(\t r)},{1.618*sin(\t r)});
\draw [shift={(1.22,0)},line width=1.5pt,color=cyan]  plot[domain=pi/3.5:0,variable=\t]({1.618*cos(\t r)},{1.618*sin(\t r)});
\node[red] at ( 0.2,  1.7)  {\tiny $2$}; 
\draw [shift={(1.22,0)},line width=0.5pt,color=cyan]  plot[domain=pi:pi/1.55,variable=\t]({-1 +1.618*cos(\t r)},{1.618*sin(\t r)});
\draw [shift={(1.22,0)},line width=1.2pt,color=red]  plot[domain=pi/1.55:pi/2.25,variable=\t]({-1+1.618*cos(\t r)},{1.618*sin(\t r)});
\draw [shift={(1.22,0)},line width=0.5pt,color=cyan, <-]  plot[domain=pi/3:pi/2.25,variable=\t]({-1 +1.618*cos(\t r)},{1.618*sin(\t r)});
\draw [shift={(1.22,0)},line width=0.5pt,color=cyan]  plot[domain=pi/3:0,variable=\t]({-1+1.618*cos(\t r)},{1.618*sin(\t r)});
\node[green] at ( 0.15,  1.45)  {\tiny $3$}; 
\draw [shift={(1.22,0)},line width=0.7pt,color=cyan, ->]  plot[domain=pi:pi/2.1,variable=\t]({-2+1.618*cos(\t r)},{1.618*sin(\t r)});
\draw [shift={(1.22,0)},line width=0.7pt,color=cyan]  plot[domain=pi/2.1:pi/2.25,variable=\t]({-2+1.618*cos(\t r)},{1.618*sin(\t r)});
\draw [shift={(1.22,0)},line width=1.2pt,color=green]  plot[domain=pi/2.25:pi/4.8,variable=\t]({-2+1.618*cos(\t r)},{1.618*sin(\t r)});
\draw [shift={(1.22,0)},line width=0.7pt,color=cyan]  plot[domain=pi/4.8:0,variable=\t]({-2+1.618*cos(\t r)},{1.618*sin(\t r)});
\node[orange] at ( -0.38,  1.05)  {\tiny $4$}; 
\draw [shift={(1.22,0)},line width=0.7pt,color=cyan, ->]  plot[domain=pi:pi/3,variable=\t]({-3+1.618*cos(\t r)},{1.618*sin(\t r)});
\draw [shift={(1.22,0)},line width=0.7pt,color=cyan]  plot[domain=pi/3:pi/4.8,variable=\t]({-3+1.618*cos(\t r)},{1.618*sin(\t r)});
\draw [shift={(1.22,0)},line width=1.2pt,color=orange]  plot[domain=pi/4.8:pi/5.3,variable=\t]({-3+1.618*cos(\t r)},{1.618*sin(\t r)});
\draw [shift={(1.22,0)},line width=0.7pt,color=cyan]  plot[domain=pi/5.3:0,variable=\t]({-3+1.618*cos(\t r)},{1.618*sin(\t r)});
\node[purple] at (  0.35,  1.02)  {\tiny $5$}; 
\draw [shift={(3.23,0)},line width=0.7pt,color=cyan, ->]  plot[domain=pi:pi/1.08,variable=\t]({2.94*cos(\t r)},{2.94*sin(\t r)});
\draw [shift={(3.23,0)},line width=0.7pt,color=cyan]  plot[domain=pi/1.08:pi/1.11,variable=\t]({2.94*cos(\t r)},{2.94*sin(\t r)});
\draw [shift={(3.23,0)},line width=1.2pt,color=purple]  plot[domain=pi/1.11:pi/1.135,variable=\t]({2.94*cos(\t r)},{2.94*sin(\t r)});
\draw [shift={(3.23,0)},line width=0.7pt,color=cyan, ->]  plot[domain=pi/1.135:pi/2.4,variable=\t]({2.94*cos(\t r)},{2.94*sin(\t r)});
%
    \begin{scope}   
        \clip (\myxlow,0) rectangle (\myxhigh,1.1);
        \foreach \i in {1,...,\myiterations}
        {   \pgfmathsetmacro{\mysecondelement}{\myxlow+1/pow(2,floor(\i/3))}
            \pgfmathsetmacro{\myradius}{pow(1/3,\i-1}
            \foreach \x in {-4,\mysecondelement,...,4}
            {   \draw[very thin] (\x,0) arc(0:180:\myradius);
                \draw[very thin] (\x,0) arc(180:0:\myradius);
            }   
        }
    \end{scope}
  \foreach \x  in { -3.5,-2.5, -1.5,-0.5, 0.5,  0.5, 1.5, 2.5, 3.5}
 {\draw (\x,0)--(\x,3.5);
 } 
    \begin{scope} 
        \begin{pgfonlayer}{background}
            \clip (-0.5,0) rectangle (0.5,3.5);
            \clip   (1,3.5) -| (-1,0) arc (180:0:1) -- cycle;
            \fill[gray,opacity=0.5] (-1,-1) rectangle (1,3.5);
        \end{pgfonlayer}
    \end{scope}
\end{tikzpicture}
\noindent
\caption{The standard tesselation of $\mathbb H$ is given by images under $\text{SL}(2, \mathbb Z)$ of its standard fundamental domain (in gray). A geodesic $\gamma$ of $\mathbb H$ and the first five geodesic arcs --- formed by intersections of  lifts of $\gamma$ with the standard fundamental domain --- which piece together to match an initial portion of the projection of $\gamma$ to $\mathcal M$. }
\label{f:standardTess}
\end{figure}

The modular surface $\mathcal M$ is the quotient of the Poincar\'e upper half-plane model $\mathbb H = \{z = x + i y\,|\, y >0\}$ of the hyperbolic plane by the group $G = SL_2(\mathbb Z)$ \footnote{Since the action of $-Id$ is trivial, one could quotient $G$ by its center $\{-Id, Id\}$ and consider the group $PSL(2,\mathbb Z)$. This changes nothing of importance, and we will ignore this detail. Similarly, we use the term surface to include possible orbifolds.} acting by fractional linear transformations.   A fundamental domain for this action is given by the triangular region bounded below by the unit circle and lying between the vertical lines of $x$-values $\pm 1/2$, choices can be made as to which portions of the boundary to include.   The images of this fundamental domain under the elements of $G$  tesselate  $\mathbb H$, let us call this the {\em standard tesselation} and these images its {\em tiles}.  First described in 1877 by Dedekind \cite{Dedekind}, the reader has almost surely seen figures representing this tesselation.   See Figure~\ref{f:standardTess}.

     The geodesics of $\mathbb H$ are the vertical half-lines and  the semi-circles with center on the real line.    Each oriented geodesic  $\gamma$ has two endpoints on the real line or at $\infty$, the {\em feet} of the geodesic,  its {\em past} foot   $\gamma_{-\infty}$ and its   {\em future} foot  $\gamma_{+\infty}$.  For ease, we use the terms geodesic to denote oriented geodesic.  The geodesics of $\mathcal M$ are the projections of the geodesics of $\mathbb H$.

Key to Series' approach is to focus not on the standard tesselation but rather upon what she dubs the {\em Farey tesselation}, see Figure~\ref{f:FareyTess}.   In briefest terms, this is the tesselation of $\mathbb H$ by the {\em Farey triangles}, the images under $G$ of the (ideal) triangle $\Delta$ of vertices $0,1$ and $\infty$.   
The edges of the Farey tesselation are exactly the images of the imaginary axis under $G$, and the vertices all lie on the line at infinity: they are exactly all rational points, and ${\infty}$.     A goal of the work is to associate to a general geodesic the infinite sequence of Farey triangles that it meets.  

There are elements $S, T \in G$ respectively acting as $z\mapsto -1/z$ and $z\mapsto z+1$.  The element $S$ maps the imaginary axis below $i$ to its portion above $i$,  the product $TS$ fixes the primitive sixth root of unity $e^{i\pi/3}$, and acts as a rotation which maps $\Delta$ to itself while permuting its edges. One finds that all of the edges of the Farey tesselation project to the projection of the ray of the imaginary axis from $i$ to $\infty$.   Let us call this $\mathcal S$.   In the paper, Series  introduces a  coding in the alphabet $\{L,R\}$ of the geodesics $\gamma$ of $\mathcal M$, where the letter reports on the geodesic segment of $\gamma$ between two consecutive intersections with $\mathcal S$: telling us whether the segment leaves with the cusp on the left or on the right, see Figure~\ref{f:FareyTess}.   This gives a unique infinite binary coding on the alphabet $\{L,R\}$ to any geodesic that does not end or start in a cusp.

The continued fractions considered in the article are the regular continued fractions, and we will also simply write continued fractions.   A continued fraction expansion is of the form  
\[ [a_0; a_1, a_2, \dots] := a_0 + \cfrac{1}{a_1 + \cfrac{1}{a_2+\cfrac{1}{\ddots}}}\;,\] 
 with $a_0 \in \mathbb Z$ and $a_n \in \mathbb N$ for all $n \in \mathbb N$.  If $x$ is an irrational real number, then $x$ can be expressed uniquely in this form,  in the sense that $x$ is the limit of its {\em convergents}, the rational numbers $p_n/q_n = [a_0; a_1, a_2, \dots, a_n]$.   The various $a_i$ are called the {\em partial quotients} of $x$.   Underlying all of this is the {\em Gauss map} $T: [0,1) \to [0,1)$ which fixes $x=0$ and otherwise sends $x$ to $1/x - \lfloor 1/x\rfloor$\,.  Identifying each $x$  with its  expansion,  the Gauss map acts as the shift map on the continued fraction expansion.

The main result  of ``{\em The modular surface and continued fractions}" is that, if the coding explained above by the Farey tesselation  is written in the form $L^{a_0}R^{a_1}L^{a_2}\ldots$,  with $a_0$ positive, then the larger  foot of the geodesic  has continued fraction expansion $ [a_0; a_1, a_2, \dots]$.     A refinement of this coding  (considering also the past of the geodesic) allows Series to determine a two dimensional   Poincar\'e section in the unit tangent bundle of the modular surface which surjects onto the unit interval such that the first return under the geodesic flow on the section  commutes with the Gauss map.   
This gives a straightforward determination of the invariant measure for the Gauss map, and indeed determines the dynamical system of the Gauss map as a factor of the first return system by geodesic flow to the explicit Poincar\'e section in the unit tangent bundle of the modular surface.   

 A very pleasant aspect of the paper is its precision --- the various sets which arise are carefully defined, and results are stated in clear topological terms.

\section{History}    
There is a rich history of the study of the connections between the geodesics of the modular surface and continued fractions.

 \subsection{Continued fractions}\label{ss:CF}    The notion of continued fractions comes directly and naturally from the Euclidean algorithm.   Given real quantities $0< a<b$, we find the number of times $n_1$ that $a$ goes into $b$ and thus determine a remainder $r_0 = b- n_1 a$ with  $0\le  r_0 < a$. If $r_0 \neq 0$, we may repeat finding $0\le a - n_2 r_0 < r_0$.   Continuing in this manner gives   
 \[ \frac{b}{a} =  n_1 + \cfrac{1}{n_2 + \cfrac{1}{\ddots}}\;,\]  and indeed one sees continued fractions emerging.

 Continued fractions themselves have been studied for centuries,  with names of some of the most celebrated  mathematicians associated to fundamental results about regular continued fractions and their generalizations, for example:  Euler, Gauss, Lagrange and Galois.

 Whereas a real number has a finite continued fraction expansion if and only if it is a rational number, any real number with a periodic expansion satisfies a quadratic polynomial with rational coefficients.  From this, periodic  continued fraction expansions correspond to real quadratic algebraic numbers.    Gauss and followers   exploited this connection to relate continued fractions to quadratic forms and were led to reduction algorithms for quadratic forms.    In particular, continued fraction expansions can be used to study class groups of real quadratic fields.   There is also a long history of studying other  fields by use of generalized continued fractions, going back to at least Minkowski's  1896 \cite{Minkowski}.

  The central nature of the Euclidean algorithm is such that there are applications of continued fractions in numerous domains, some apparently far afield from the geodesic flow on the modular surface.    Examples include:  Conway's 1970 \cite{Conway} introduction of their use in knot theory; and,  Hirzebruch's  1953 \cite{HirzebruchQuotient} use  of finite expansions in terms of a variant of continued fractions in resolving quotient singularities of algebraic surfaces, and his 1971 \cite{HirzebruchCusp} application of  periodic expansions in the resolution of cusp singularities.

   Continued fractions are also inherently related to Diophantine approximation, with (regular) continued fraction convergents giving the best rational numbers approximating a nonzero real number $x$ in the sense that if $p/q$ is a rational number with $|x - p/q| < 1/(2 q^2)$ then $p/q = p_n/q_n$ is the $n^{\text{th}}$ convergent of $x$ for some $n$.       Due to this, restrictions on continued fraction expansions of a number are part of the hypotheses in various theorems  throughout mathematics.

   The dynamical  view of continued fractions also has a long tradition, arguably beginning with Gauss'   determination of the invariant measure   of the Gauss map, $d\mu=\frac 1{\ln 2}\frac {1}{x+1}\, dx$. Intriguingly,  it remains a mystery as to how Gauss discovered the invariant measure.   The system of the Gauss map has long been viewed as a fundamental example in dynamics and ergodic theory.

\subsection{Coding of geodesics} \label{ss:coding} 

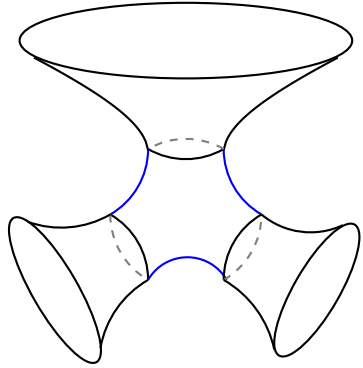
\begin{figure}[h!]
    \centering
    \begin{tikzpicture}[scale = 1]  
        \draw[thick, blue] (1,0) arc (240:180:1);
        \draw[thick] (1/2,{sqrt(3)/2}) arc (300:240:1);
        \draw[thick, blue] (-1/2,{sqrt(3)/2}) arc (0:-60:1);
        \draw[thick] (-1,0) arc (60:0:1);
        \draw[thick, blue] (-1/2,{-sqrt(3)/2}) arc (150:30:0.6);
        \draw[thick] (1/2,{-sqrt(3)/2}) arc (180:120:1);
          \draw[domain=2.24:3.75,samples=50,thick] plot(\x-1.74, {sqrt(3)/2 +0.95*sqrt(\x^2/5-1)-0.07});
          \draw[domain=2.24:3.75,samples=50,thick] plot(-\x+1.74, {sqrt(3)/2 +0.95*sqrt(\x^2/5-1)-0.07});
         %
        %
        \draw[thick, gray, dashed] (1/2,{sqrt(3)/2}) arc (60:120:1);
        \draw[thick, gray, dashed] (-1,0) arc (180:240:1) node[pos=0.5, anchor=east] { };
        \draw[thick, gray, dashed] (1/2,{-sqrt(3)/2}) arc (300:360:1) node[pos=0.7, anchor=north, xshift=2mm] { }; 
     \draw[thick] (1/2,{-sqrt(3)/2}) arc (60:12:1.4);
     \draw[thick] (1,0) arc (230:295:1.01);
      \draw[thick] (-1/2,{-sqrt(3)/2}) arc (120:170:1.3);
      \draw[thick] (-1,0) arc (300:250:1.3);
      \draw[thick, rotate=90] (2.3,0) ellipse (0.5 and 2.2);
      \draw[thick,  rotate=210] (2,0) ellipse (0.3 and 1.1);
      \draw[thick,  rotate=330] (2,0) ellipse (0.3 and 1.0);
    \end{tikzpicture}  
    \caption{  A pair of pants with three funnels and three `seams' (in blue) along a central compact region.  Geodesics exiting the central compact region remain in a funnel.}
 \label{f:PantsWithSeams}
\end{figure}

 Following Ghys \cite{Ghys}, we indicate some of Hadamard's insights in his  aforementioned work of 1898.    The simplest type of  surfaces which Hadamard considered,  now called `pairs of pants',  can be visualized as deleting three disks from  a sphere  and glueing to each boundary circle an infinite cylinder which flairs out so as to have infinite hyperbolic area, see Figure~\ref{f:PantsWithSeams}.   Hadamard showed that each of these `funnels' can be separated from the remainder of the surface by a unique geodesic and that any geodesic which crosses this separating geodesic wanders off into the funnel never to return.  Hadamard considered the set of closed geodesics which remained in the finite portion of the surface.  To this end he cut the surface along a further set of geodesics (in the pair of pants setting, he took  three separating geodesic segments which run along `seams' from one funnel to another), and assigned to a geodesic the bi-infinite sequence of labels of the separating curves which the geodesic crossed.  Besides proving the density of the closed geodesics among the bounded geodesics, Hadamard also showed  that at each point of the finite portion of the surface the set of unit tangent vectors to bounded geodesics passing through the given point form a Cantor set.   
 See  the  excellent (mainly) expository article of Ghys for more on Hadamard's work.
 
        In 1927,   Koebe \cite{Koebe} used true cutting sequences.  To a general geodesic he associated a bi-infinite sequence of elements of the surface's fundamental group, where each element corresponded to the projection of a side of a fixed fundamental region for the surface's fundamental group acting on the hyperbolic plane.  (In fact, Koebe writes that in 1917 he had notes for this procedure.)   In works culminating in his 1938 notes \cite{MorseNotes}, Morse also came to the use of cutting sequences.   (Series  \cite{Series86} explains why  Morse's earlier \cite{Morse1} and \cite{Morse2}  are not applications of true cutting sequences.) 
 
 Explicitly relating the  geodesics  of the modular surface to continued fractions was first done by E.~Artin \cite{Artin} in 1924. Artin shows the existence of geodesics each of whose set of unit tangent vectors  is  dense in the {\em unit tangent bundle}, thus in the collection of all unit tangent vectors, of $\mathcal M$.   Artin writes that the motivation for the work  arose from correspondence with his former Ph.D. thesis advisor,   Herglotz.   The only background literature mentioned in the paper is a textbook on differential geometry by  Blaschke, so we cannot know if Artin was aware of the work  of Hadamard.     Recall that each geodesic of $\mathcal M$  is the projection of a semi-circle (or vertical line) whose feet   lie on the extended real axis.  Indeed, each geodesic has infinitely many of these {\em lifts}.   To each lift, Artin associated the continued fraction expansions of its feet.  Restricting to a special subset of the geodesics, he was able to show the desired denseness result. 
 
 See works of S.~Katok, both singly and with Ugarcovici, in particular \cite{KUbAMS}, for excellent discussion of coding of geodesics. 
They 
refer  to the two types of coding exemplified by the  approaches of Koebe and Morse and respectively by Artin as {\em geometric} and {\em arithmetic}.       In  geometric coding  one marks a fundamental domain and  assigns to a geodesic the bi-infinite sequence of group elements identifying the outgoing side  for a lift of the geodesic with the incoming side for the subsequent lift.   In arithmetic coding one associates to a geodesic those elements of the group  which successively bring  a lift of the geodesic into a predetermined  `reduced set'.  Also  called `reduction algorithms' \cite{KUbAMS}, this arithmetic approach has the longer history.  Indeed, this is directly related to Gauss' reduction of quadratic forms.    
 
 Let us say a bit more about the approach of Artin and its variants.  Again due to the aforementioned elements $S, T \in G$,  the projection to $\mathcal M$ of the edges of the standard   tesselation is the union of two curves: $a$,   the projection  of the vertical line joining  $e^{i\pi/3}$  to $\infty$;   and $b$,  the projection of  the arc joining   $i$ to $e^{i\pi/3}$.  Let us call the projections of  $i$ and $e^{i\pi/3}$ simply `special points'.     To a  geodesic   of $\mathcal M$, one  associates the  cutting sequence given by the sequence of $a$ and $b$ cut by the geodesic.  An obvious difficulty arises  if the geodesic hits either of the special points.  Ignoring this difficulty, one could create a corresponding symbolic dynamical system given by the  {\em vertex} shift   on the graph whose vertices are labeled $a$ and $b$ and   whose  edges join  these.     The exact determination of the rules identifying the edges in the graph is a delicate matter.  
 
 In \cite{Series}, Series points to her 1981 work \cite{Series81} and to Adler-Flatto's 1982 \cite{AdlerFlatto} as giving ways to resolve   the  difficulties remaining in Artin's approach, but these approaches are quite cumbersome, and not as natural as the approach of \cite{Series}.

The coding of geodesics lies at the very foundation of dynamical systems and ergodic theory.     In 1935 Hedlund \cite{HedlundMod} used Artin's approach to prove the ergodicity of the geodesic flow on the unit tangent bundle of $\mathcal M$.  In \cite{HedlundClosedSurfs} he proved the analogous result for the closed surfaces considered by Nielsen, using Nielsen's boundary expansions. (In short,  boundary expansions are forward orbits of points of the boundary of the hyperbolic plane under a function induced by the fractional linear action of elements pairing the sides of a chosen fundamental domain for the Fuchsian group uniformizing the surface.)   Soon thereafter,  E.~Hopf proved the ergodicity of the geodesic flow on the unit tangent bundle of any finite area hyperbolic surface;   he revisits this, with added detail,  in  \cite{Hopf}.   (See  Coven and Nitecki \cite{CovenNitecki}  for an iconoclastic discussion of the early history of symbolic dynamics.)

\subsection{Farey tesselation}\label{ss:FareyTess}  The interpretation of continued fractions in terms of the Farey tesselation goes back to H.~J.~Smith in 1877, but seems to have lain fallow until it was greatly furthered by  Humbert in 1916.    (Note that while  \cite{Series} cites short articles in 1916 of Humbert,  the longer article \cite{Humbert} of the same year and title(!)   more fully addresses matters.)    Humbert relates the  continued fraction expansion of a real number $x$ to the sequence of intersections of the infinite vertical ray with endpoint at $x$ with the Farey triangle tesselation.    Whenever this vertical ray intersects a Farey triangle (other than the triangle with vertex at infinity), it isolates one of the vertices from the others.   Humbert calls this a {\em tip} related to $x$, each such is a rational number, and shows that the continued fraction convergents of $x$ are exactly the tips so formed.   The word geodesic is not to be found in that article.  Still, Humbert studies those semi-circles   whose feet are a real quadratic surd  and its conjugate.    The projections of these semi-circles are exactly the closed geodesics on the modular surface.    As for the vertical ray, the tips determined by the intersection of  such a semi-circle with the Farey tesselation gives continued fraction convergents.  Furthermore,  Humbert calls the `order' of a convergent of $x$ the number of Farey triangles for which it is the tip.  He shows that the order is, in the appropriate sense,  the  next partial quotient and uses this to recover results of Lagrange and Galois about periodicity of continued fraction expansions,  as well as results related to class numbers of indefinite  integral binary forms.    
 
 Curiously enough,  the Farey tesselation  is most naturally associated to a degree three cover of the modular surface.  Indeed, the aforementioned ideal triangle $\Delta$ is a fundamental domain for the subgroup $G^3$ generated by the cubes of the elements of $G$.  (This subgroup is one of only two finite index normal subgroups of $G$ which is not a free group.)   
The quotient surface is a punctured sphere with three order two quotient singularities.   This has a double cover which is the punctured torus that is the quotient of $\mathbb H$ by the commutator subgroup of $G$, and it also has a degree four cover by the  four times punctured sphere which is the quotient of $\mathbb H$ by the congruence subgroup modulo three.    These three surfaces were shown in the 1980s to share sets of simple closed geodesics  (in the sense that the lift to $\mathbb H$  of  a simple closed geodesic on any of these surfaces projects to each of the three surfaces to be simple and closed), with the importance of this being that the Markoff spectrum of Diophantine approximation is naturally in one-to-one correspondence with the elements of any of these sets.  See \cite{SeriesGeomMarkoffNumbers} for a discussion.

Series writes in \cite{Series} that she was introduced to the Farey tesselation by Moeckel's 1982  \cite{Moeckel}.   There, certain number theoretic questions are answered by considering what Moeckel called {\em excursions into the cusp} of general geodesics ---  in Humbert's terms each such excursion is the geodesic segment corresponding to a single tip (thus convergent), and each is partitioned into subsegments by the Farey triangles having this tip as vertex.   That is, the geodesic coding aspect of Moeckel's work is that of Humbert, except that Moeckel  considers also non-closed geodesics.\\   

 From this brief history, we see that, by the beginning of the 1980s, it had been long recognized that a strong connection linked the geodesic flow on the modular surface and the regular continued fractions, but a satisfying explicit formulation of this link was still missing.  
 
 \section{The three main theorems}  
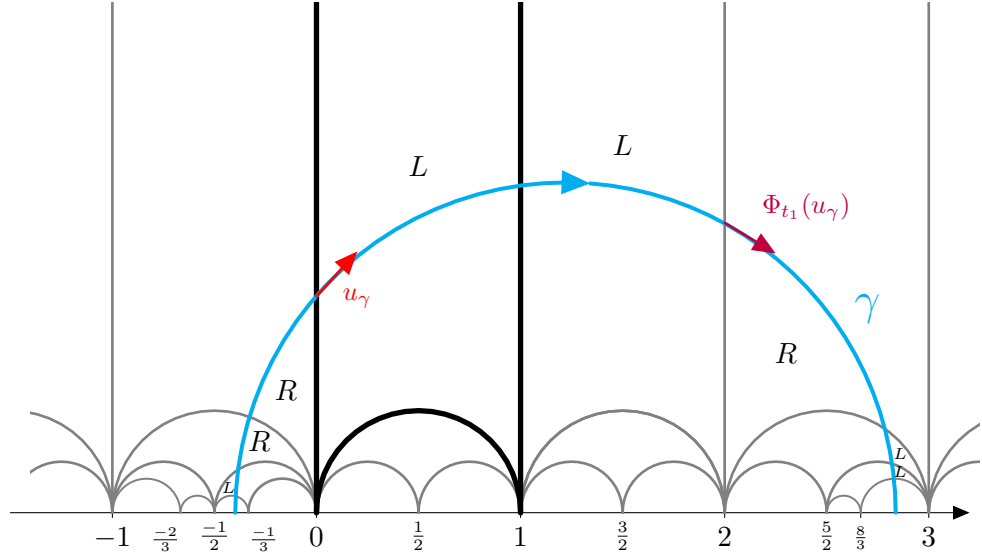
\begin{figure}[h]  
\begin{tikzpicture}[line cap=round,line join=round,>=triangle 45,x=3.0cm,y=3.0cm, scale=0.9]
\def\heightConstant{2.5}
\draw[->,color=black] (-1.5,0) -- (3.2,0);
\foreach \x in {-1,0,1,2,3}
\draw[shift={(\x,0)},color=black] (0pt,2pt) -- (0pt,-2pt) node[below] {\large $\x$};
\foreach \x in {-1, 1, 3, 5}
\draw[shift={(\x/2,0)},color=black]  node[below] {\small $\frac{\x}{2}$}; 
\node [label={[shift={(-2/3, -0.25)}]\tiny $\frac{-2}{3}$}] {};
\node [label={[shift={(-1/3+0.1, -0.25)}]\tiny $\frac{-1}{3}$}] {};
\draw[shift={(8/3,0)},color=black] (0pt,2pt) -- (0pt,-2pt) node[below] {\tiny $\frac{8}{3}$};
\clip(-1.4,-0.1) rectangle (3.25,\heightConstant);
%
\draw [shift={(1.5,0)},line width=1pt,color=gray]  plot[domain=0:pi,variable=\t]({ 0.5*cos(\t r) },{0.5*sin(\t r)});
\draw [shift={(0.75,0)},line width=1pt,color=gray]  plot[domain=0:pi,variable=\t]({0.25*cos(\t r)},{0.25*sin(\t r)});
\draw [shift={(0.25,0)},line width=1pt,color=gray]  plot[domain=0:pi,variable=\t]({0.25*cos(\t r)},{0.25*sin(\t r)});
\draw [shift={(1.5,0)},line width=1pt,color=gray]  plot[domain=0:pi,variable=\t]({-0.5*cos(\t r)},{0.5*sin(\t r)});
\draw [shift={(1.25,0)},line width=1pt,color=gray]  plot[domain=0:pi,variable=\t]({-0.25*cos(\t r)},{0.25*sin(\t r)});
\draw [shift={(1.75,0)},line width=1pt,color=gray]  plot[domain=0:pi,variable=\t]({0.25*cos(\t r)},{0.25*sin(\t r)});
\draw [shift={(-0.5,0)},line width=1pt,color=gray]  plot[domain=0:pi,variable=\t]({-0.5*cos(\t r)},{0.5*sin(\t r)});
\draw [shift={(-0.75,0)},line width=1pt,color=gray]  plot[domain=0:pi,variable=\t]({-0.25*cos(\t r)},{0.25*sin(\t r)});
\draw [shift={(-1.5,0)},line width=1pt,color=gray]  plot[domain=0:pi,variable=\t]({0.5*cos(\t r)},{0.5*sin(\t r)});
\draw [shift={(-1.25,0)},line width=1pt,color=gray]  plot[domain=0:pi,variable=\t]({0.25*cos(\t r)},{0.25*sin(\t r)});
\draw [shift={(-1.75,0)},line width=1pt,color=gray]  plot[domain=0:pi,variable=\t]({-0.25*cos(\t r)},{0.25*sin(\t r)});
\draw [shift={(-0.25,0)},line width=1pt,color=gray]  plot[domain=0:pi,variable=\t]({0.25*cos(\t r)},{0.25*sin(\t r)});
\draw [shift={(-0.167,0)},line width=1pt,color=gray]  plot[domain=0:pi,variable=\t]({0.167*cos(\t r)},{0.167*sin(\t r)});
\draw [shift={(2.25,0)},line width=1pt,color=gray]  plot[domain=0:pi,variable=\t]({-1*0.25*cos(\t r)+0*0.25*sin(\t r)},{0*0.25*cos(\t r)+1*0.25*sin(\t r)});
\draw [shift={(2.75,0)},line width=1pt,color=gray]  plot[domain=0:pi,variable=\t]({1*0.25*cos(\t r)+0*0.25*sin(\t r)},{0*0.25*cos(\t r)+1*0.25*sin(\t r)});
\draw [shift={(3.25,0)},line width=1pt,color=gray]  plot[domain=0:pi,variable=\t]({-1*0.25*cos(\t r)+0*0.25*sin(\t r)},{0*0.25*cos(\t r)+1*0.25*sin(\t r)});
\draw [shift={(3.75,0)},line width=1pt,color=gray]  plot[domain=0:pi,variable=\t]({1*0.25*cos(\t r)+0*0.25*sin(\t r)},{0*0.25*cos(\t r)+1*0.25*sin(\t r)});
\draw [shift={(3.5,0)},line width=1pt,color=gray]  plot[domain=0:pi,variable=\t]({-1*0.5*cos(\t r)+0*0.5*sin(\t r)},{0*0.5*cos(\t r)+1*0.5*sin(\t r)});
\draw [shift={(2.5,0)},line width=1pt,color=gray]  plot[domain=0:pi,variable=\t]({-1*0.5*cos(\t r)+0*0.5*sin(\t r)},{0*0.5*cos(\t r)+1*0.5*sin(\t r)});
\draw [shift={(-5/12,0)},line width=0.75pt,color=gray]  plot[domain=0:pi,variable=\t]({(1/12)*cos(\t r)},{(1/12)*sin(\t r)});
\draw [shift={(-5/6,0)},line width=0.75pt,color=gray]  plot[domain=0:pi,variable=\t]({(1/6)*cos(\t r)},{(1/6)*sin(\t r)});
\draw [shift={(17/6,0)},line width=0.75pt,color=gray]  plot[domain=0:pi,variable=\t]({(1/6)*cos(\t r)},{(1/6)*sin(\t r)});
\draw [shift={(-7/12,0)},line width=0.75pt,color=gray]  plot[domain=0:pi,variable=\t]({(1/12)*cos(\t r)},{(1/12)*sin(\t r)});
\draw [shift={(31/12,0)},line width=0.75pt,color=gray]  plot[domain=0:pi,variable=\t]({(1/12)*cos(\t r)},{(1/12)*sin(\t r)});
%
\draw [line width=2pt,color=black] (0,0) -- (0,\heightConstant);
\draw [line width=2pt,color=black] (1,0) -- (1,\heightConstant);
\draw[shift={(0.5,0)},line width=2pt,color=black]  plot[domain=0:pi,variable=\t]({1*0.5*cos(\t r)+0*0.5*sin(\t r)},{0*0.5*cos(\t r)+1*0.5*sin(\t r)});
\draw [line width=1pt,color=gray] (-1,0) -- (-1,\heightConstant);
\draw [line width=1pt,color=gray] (-3,0) -- (-3,\heightConstant);
\draw [line width=1pt,color=gray] (3,0) -- (3,\heightConstant);
\draw [line width=1pt,color=gray] (-3,0) -- (-3,\heightConstant);
\draw [line width=1pt,color=gray] (-2,0) -- (-2,\heightConstant);
\draw [line width=1pt,color=gray] (2,0) -- (2,\heightConstant);
\draw [line width=1pt,color=gray] (4,0) -- (4,\heightConstant);
\draw [shift={(1.22,0)},line width=1.5pt,color=cyan, ->]  plot[domain=pi:pi/2.1,variable=\t]({1.618*cos(\t r)},{1.618*sin(\t r)});
\draw [shift={(1.22,0)},line width=1.5pt,color=cyan]  plot[domain=pi/2.1:0,variable=\t]({1.618*cos(\t r)},{1.618*sin(\t r)});
\node[cyan] at (2.7, 1)  {\huge $\gamma$};
\node at (0.5, 1.7)  {\large $L$};
\node at (1.5, 1.8)  {\large $L$};
\node at (2.3, 0.78)  {\large $R$};
\node at (2.86, 0.29)  {\tiny $L$};
\node at (2.86, 0.2)  {\tiny $L$};
\node at (-0.43, 0.125)  {\tiny $L$};
\node at (-0.28, 0.35)  {\large $R$};
\node at (-0.15, 0.6)  {\large $R$};
\draw[red, line width=1pt,->] (0,1.06)--(0.2, 1.28) {};
\node[red] at (0.2, 1.05) {$u_{\gamma}$};
%
\draw[purple, line width=1pt,->] (2,1.42)--(2.25, 1.27) {};
\node[purple] at (2.4, 1.5) {$\Phi_{t_1}(u_{\gamma})$};
\end{tikzpicture}
\noindent
\caption{The ideal triangle (edges in black),  further triangles in the Farey tesselation (edges in gray), a geodesic ($\gamma$ in cyan) and a portion of its $L, R$-sequence,   the unit tangent vector  $u_{\gamma}$ (in red) of base point the intersection of $\gamma$ with the imaginary axis,  and the image $\Phi_{t_1}(u_{\gamma})$ of this tangent vector under the geodesic flow for a time $t_1$. Note that $\gamma$ here is the same geodesic as $\gamma$ in Figure~\ref{f:standardTess}.} 
\label{f:FareyTess}
\end{figure}

 Series  \cite{Series}   perfected the  use of the Farey tesselation by leveraging the fact that the tips of  consecutive excursions into cusps  are alternately to the right or to the left as determined by the oriented geodesic under consideration.   To each non-vertical geodesic $\gamma$ of $\mathbb H$ (whose feet are not rational numbers) she associated a bi-infinite sequence of the letters $R$ and $L$, determined by whether the corresponding tip lies to the right or to the left, see Figure~\ref{f:FareyTess}.   (If either foot is a rational number, then the sequence correspondingly begins or ends with a special symbol.)    Thus, Series' key insight was to recognize that Moeckel's  excursions into the cusp naturally came in two types,   completely changing how we now read this figure.   Compare  [\cite{Moeckel}, Figure~6] and [\cite{SeriesArtInBook}, Figure~7].
 
  Collecting consecutive letters  of the same type (called {\em runs}), results in a bi-infinite sequence of natural numbers.   A geodesic {\em changes type} at a point on the edge of a Farey triangle if its associated sequence changes letter there.      Since each geodesic of $\mathcal M$ can be lifted to $\mathbb H$, the above gives     a bi-infinite sequence in $\{L,R\}$ determined by behavior between intersections with $\mathcal S$ mentioned in the introduction above.   Whereas the Artin approach gives  a vertex shift symbolic system, here the coding can be seen as a {\em edge shift} symbolic system:    the curve $\mathcal S$ corresponds to a single vertex and there are two edges labeled as $R, L$.

  The paper \cite{Series} has three main results, listed as Theorems A through C.  These  results combine to show, as the Zentralblatt reviewer of \cite{Series} wrote:   ``The geodesic flow is a special flow over a shift which is essentially the natural extension of the continued fraction transformation, under a height function given by crossing times across   [tiles of the Farey tesselation]." 
 \bigskip

  \subsection{Decorated sequences and a projection, Theorem~A} 
      Let $A$ be the set of all  geodesics of $\mathbb H$ whose feet are of opposite signs, with $|\gamma_{+\infty}| \ge 1$ and $0<|\gamma_{-\infty}| \le 1$.   Each geodesic $\gamma \in A$   meets the imaginary axis so as to change type there.  Let $u_{\gamma}$ be the unit vector tangent to $\gamma$ at this point of intersection, see Figure~\ref{f:FareyTess}.  Consider the function sending $\gamma$ to the unit tangent vector of $\mathcal M$ which is the projection of $u_{\gamma}$.        The projection belongs to the set $X$ of unit tangent vectors $v$ with basepoint $x$ on $\mathcal S$   and such that the lift to $\mathbb H$ of the geodesic that $v$ defines changes type at (the lift of) $x$.     Theorem A states that this function  from $A$ to $X$ is  surjective, continuous and open with respect to the natural topologies. 
  
 Furthermore,   to each  $\gamma \in A$     \cite{Series}  associates a (usually doubly infinite) sequence of natural numbers $\dots, n_{-1}, n_0, n_1, n_2, \dots$ recording the length of runs of each of the letters $R, L$ in the sequence for $\gamma$ and with $n_1$ corresponding to the run associated to the segment of $\gamma$ directly after the basepoint of $u_{\gamma}$.  In order to not lose the information of which letter has its run of length $n_1$,  later one associates to $\gamma$ an ordered pair whose first entry is the sequence as above and  whose second entry is a $0$ if that run is of the letter $L$ and a $1$ otherwise.   Theorem~A also states that  when that determining letter is $L$,  then  $\gamma_{+\infty}$ has continued fraction expansion  $[n_1, n_2, \dots]$ while $-1/\gamma_{-\infty}$ has continued fraction expansion $[n_0, n_{-1}, \dots]$.  A similar statement is given for the case when $L$ is replaced by $R$.

\begin{figure}[h]
\scalebox{0.68}{
\noindent
\begin{tabular}{lcr}
\begin{tikzpicture}[line cap=round,line join=round,>=triangle 45,x=3.0cm,y=3.0cm, scale=0.9]
\def\heightConstant{2.5}
\draw[color=black] (-1.2,0) -- (1.2,0);
\draw [line width=2pt,color=black] (0,0) -- (0,\heightConstant);
\draw [line width=1pt,color=black] (1,0) -- (1,\heightConstant);
\draw[shift={(0.5,0)},line width=1pt,color=black]  plot[domain=0:pi,variable=\t]({1*0.5*cos(\t r)+0*0.5*sin(\t r)},{0*0.5*cos(\t r)+1*0.5*sin(\t r)});
\draw [line width=1pt,color=black] (-1,0) -- (-1,\heightConstant);
\draw [line width=1pt,color=black] (0,0) -- (0,\heightConstant);
\draw[shift={(-0.5,0)},line width=1pt,color=black]  plot[domain=0:pi,variable=\t]({1*0.5*cos(\t r)+0*0.5*sin(\t r)},{0*0.5*cos(\t r)+1*0.5*sin(\t r)});
\draw[shift={(-1,0)},color=black] (0pt,2pt) -- (0pt,-2pt) node[below] {\large $n_1-1$};
\draw[shift={(0,0)},color=black] (0pt,2pt) -- (0pt,-2pt) node[below] {\large $n_1$};
\draw[shift={(1,0)},color=black] (0pt,2pt) -- (0pt,-2pt) node[below] {\large $n_1+1$};
\draw [shift={(1.22-2,0)},line width=1.5pt,color=cyan, ->]  plot[domain=pi/1.5:pi/8,variable=\t]({1.618*cos(\t r)},{1.618*sin(\t r)});
\draw [shift={(1.22-2,0)},line width=1.5pt,color=cyan]  plot[domain=pi/8:0,variable=\t]({1.618*cos(\t r)},{1.618*sin(\t r)});
\node[cyan] at (0.7, 1.1)  {\huge $\gamma$};
\node[black]  at (-0.5, 1.8)  {\large $L$};
\node[black] at (0.29, 0.78)  {\large $R$};
\draw[purple, line width=1pt,->] (0,1.42)--(0.25, 1.27) {};
\node[purple] at (0.2, 1.5) {$v$};
\end{tikzpicture}
&
\begin{tikzpicture}[x=1.5cm,y=5cm] 
\node at (0, 0) {\phantom{here}};
\draw[->, ultra thick, black] (0.2, 0.5)--(0.8, 0.5) -- (.7, 0.55) -- (1.3,0.55);
\node[black] at (.7, .7) {{$z\mapsto M\cdot z$}};
\node[black] at (.7, .3) {{$M=  \begin{pmatrix} 0&-1\\1&-n_1\end{pmatrix}$}};
\end{tikzpicture}
&
\begin{tikzpicture}[line cap=round,line join=round,>=triangle 45,x=3.0cm,y=3.0cm, scale=0.9]
\def\heightConstant{2.5}
\draw[color=black] (-2,0) -- (1.2,0);
\draw [line width=2pt,color=black] (0,0) -- (0,\heightConstant);
\draw [line width=1pt,color=black] (1,0) -- (1,\heightConstant);
\draw[shift={(0.5,0)},line width=1pt,color=black]  plot[domain=0:pi,variable=\t]({1*0.5*cos(\t r)+0*0.5*sin(\t r)},{0*0.5*cos(\t r)+1*0.5*sin(\t r)});
\draw [line width=1pt,color=black] (-1,0) -- (-1,\heightConstant);
\draw [line width=1pt,color=black] (0,0) -- (0,\heightConstant);
\draw[shift={(-0.5,0)},line width=1pt,color=black]  plot[domain=0:pi,variable=\t]({1*0.5*cos(\t r)+0*0.5*sin(\t r)},{0*0.5*cos(\t r)+1*0.5*sin(\t r)});
\draw[shift={(-1,0)},color=black] (0pt,2pt) -- (0pt,-2pt) node[below] {\large $-1$};
\draw[shift={(0,0)},color=black] (0pt,2pt) -- (0pt,-2pt) node[below] {\large $0$};
\draw[shift={(1,0)},color=black] (0pt,2pt) -- (0pt,-2pt) node[below] {\large $1$};
\draw [shift={(-0.63,0)},line width=1.5pt,color=cyan]  plot[domain=pi:pi/2,variable=\t]({1.05*cos(\t r)},{1.05*sin(\t r)});
\draw [shift={(-0.63,0)},line width=1.5pt,color=cyan, <-]  plot[domain=pi/2:pi/10,variable=\t]({1.05*cos(\t r)},{1.05*sin(\t r)});
\node[cyan] at (-1.5, 1.1)  {\huge $M\cdot\gamma$};
\node[black]  at (-0.5, 1.8)  {\large $R$};
\node[black] at (0.1, 0.5)  {\large $L$};
\draw[red, line width=1pt,->] (0,0.85)--(-0.2, 1) {};
\node[red] at (-0.25, 0.8) {$M\cdot v$};
\end{tikzpicture}
\end{tabular}
}
\caption{The map $z \mapsto -1/(z - n_1) = \begin{pmatrix} 0&-1\\1&-n_1\end{pmatrix}\cdot z$ sends a geodesic $\gamma$ with $\gamma_{+\infty} = [n_1, n_2,  \dots ]$ to a geodesic  whose future foot has continued fraction expansion  $-[n_2, \dots ]$.  The map sends  the Farey tile of vertices $(n_1, n_1-1, \infty)$ to the standard Farey triangle of vertices $(\infty, 1, 0)$, and the tile of vertices $(n_1, n_1+1, \infty)$ to the tile of vertices $(\infty, -1, 0)$.  It sends the unit tangent vector $v$ to $\gamma$ with basepoint on the vertical line at $x = n_1$ to a unit tangent vector with basepoint on the imaginary axis.    Compare this  to Figure~\ref{f:FareyTess} while setting $v = \Phi_{t_1}(u_{\gamma})$. 
}%
\label{f:oneMove}%
\end{figure}
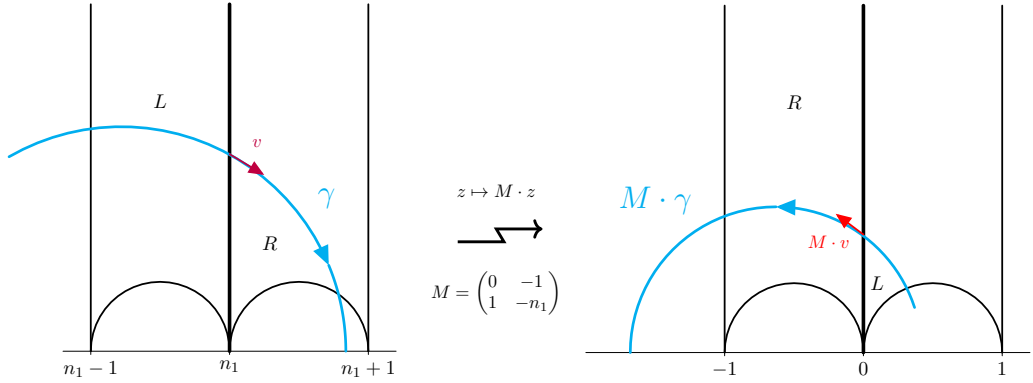
 
\subsection{Geodesic flow and a coding of a first return map, Theorem~B}  
 Theorem  B  addresses the dynamics of the situation, and in particular involves the geodesic flow on  the unit tangent bundle of $\mathcal M$.  The {\em geodesic flow} on the unit tangent bundle of $\mathbb H$ is a group action of the additive group of real numbers.  Given a real $t\ge 0$ and a unit tangent vector $u$,  since $u$ uniquely determines a geodesic $\gamma$ and a positive direction,  beginning at the basepoint of $u$ we can measure off a geodesic segment of $\gamma$  of arclength $t$ in this positive direction. We let $\Phi_t(u)$ be the unit tangent vector whose basepoint is the second endpoint of this geodesic segment.  When $t<0$ we use a segment in the opposite direction along $\gamma$.    The  maps $(t,u) \mapsto \Phi_t(u)$ give the geodesic flow.   All of this descends so as to give a well-defined geodesic flow on the unit tangent bundle of $\mathcal M$.    
   
   Let $X^*$ be the subset of those $v \in X$ such that the geodesic defined by $v$ does not exit into the cusp before returning to $\mathcal S$.  The {\em first return map} $P: X^* \to X$ sends  $v\in X^*$ to its first image  under the geodesic flow  (thus given by least positive  $t$) which also lies in $X$.  Theorem~B gives a symbolic coding of the return map $P$.

      For this, Series introduces the set $\Sigma$ of all pairs consisting of a sequence of natural numbers with indices from $N_1$ to $N_2$ where $-\infty \le N_1 \le 0 < N_2 \le \infty$ and an element of $\mathbb Z/2$.  She also considers the  subset $\Sigma^*$,  consisting of those pairs whose sequence has $N_2>1$.   She defines a self-map $\hat{\sigma}$ of $\Sigma$ which is the standard right shift on the sequences and addition modulo $2$ on the second entry.  She then defines a map $D: \Sigma \to X$,  which restricts to send $\Sigma^*$ to $X^*$ and such that  $P\circ D = \hat{\sigma} \circ D$.   Theorem~B precisely gives the topological properties   of this symbolic coding of the first return map to $X$ under the geodesic flow. In particular, with respect to natural topologies, $D$ is a continuous surjection.
   
 We say a bit about how $D$ comes about.    Given $\gamma \in A$, there is  a first time $t_1>0$ such that $\Phi_{t_1}(u_{\gamma})$ has its basepoint on the boundary of a tile of the Farey tesselation and such that $\gamma$ changes type there, again see Figure~\ref{f:FareyTess}.   For our ease, let us denote the segment of $\gamma$ from the basepoint of $u_{\gamma}$ to the basepoint of $\Phi_{t_1}(u_{\gamma})$ by $\ell$.    There is an element of $G$ sending the pair of geodesic $\gamma$ and its unit tangent vector $\Phi_{t_1}(u_{\gamma})$ to a pair of $\gamma'$ and  $u' = u_{\gamma'}$ such that $u'$ has basepoint on the imaginary axis and the sequence of natural numbers associated to $\gamma' \in A$ is simply the right shift applied to the sequence associated to $\gamma$.   See Figure~\ref{f:oneMove}    for a hint of this.  Since $\gamma'$ is a $G-$image of $\gamma$, these geodesics project to the same geodesic on $\mathcal M$.  Indeed,   the forward continuation on $\mathcal M$  of the projection of the segment  $\ell$ is the projection of the ray along $\gamma'$ beginning at the basepoint of $u'$.     From this, each unit tangent vector for $\mathcal M$ in the set $X^*$ lifts to some $u_{\gamma}$ having an associated element of $\Sigma^*$.   The first return map $P$ gives an element of $X$ that is the projection of $\Phi_{t_1}(u_{\gamma})$ and thus which is the projection of $u'$.  Of course the map $D$ goes in the opposite direction, and there are details to carefully argue.

\subsection{The Gauss map as a factor, Theorem~C}   
   Theorem~C gives an explicit ``projection" from $\Sigma$  to  $[0,1)$  with which the  
 shift on $\Sigma$ and the Gauss map on  $[0,1)$ form a commutative diagram.     
 The projection is given by  sending  a general element $(\,(n_i)_{N_1\le i\le  N_2}, w\,)$ of $\Sigma$ to  $1/[n_1, n_2, \dots]$; the notation in  \cite{Series} for this projection is  $J\circ p^+$.    Since the Gauss map $T$ sends  $1/[n_1, n_2, \dots]$ to $1/[ n_2, \dots]$  one indeed finds that $T \circ (J\circ p^+) = (J\circ p^+)\circ \hat{\sigma}$.

 \section{  Alternate codings  of the geodesic flow on the modular surface}\label{s:VarCodings}       
  We mention some of the approaches to relating geodesics of $\mathcal M$  and continued fractions  {\em not} taken in \cite{Series}.   One approach was introduced by L.~R.~Ford  in 1917.    He  used the images under elements of $G$ of the `horocycle' of $\mathbb H$ given by $y = 1$ to give a straightforward geometric argument determining the so-called Hurwitz constant for Diophantine approximation by continued fractions.  An element $M \in G$ sends that horocycle  to a Euclidean circle tangent at $p/q = M\cdot \infty$ of radius $1/(2 q^2)$.  Thus the vertical geodesic based at a real $x$ meets this circle exactly when $|x - p/q|\le 1/(2 q^2)$. Such a circle is now called a {\em Ford circle}.       Using this, one can study the continued fraction expansion of $x$; Ford's 1938  \cite{Ford1938} gives a wonderful introduction to these ideas.    However,  Ford never quite discusses codings of geodesics of $\mathcal M$.

    Adler and Flatto followed up their aforementioned 1982 paper with several more works on related situations.   Throughout,  they take a variant of the geometric approach by fixing a fundamental domain and taking as their section of the unit tangent bundle of $\mathcal M$ (or the appropriate surface) the projection of the unit tangent vectors of $\mathbb H$ pointing outwards from a chosen domain.     In the setting of $\mathcal M$, where  Series' chose to take {\em all} (non-vertical) unit tangent vectors along the geodesic given by the projection of the imaginary axis, Adler and Flatto's choice to select a subset of unit tangent vectors depending on the location of the basepoint along their boundary causes them to work quite hard to obtain a nice projection map to the unit interval.   That said, by carefully linearizing the first return map to the cross section they are able to obtain an underlying symbolic sequence that is a Markov chain.

    In 1994 Arnoux \cite{ArnouxCodage}  specialized ideas of W.~Veech for geometrically representing the so-called Teichm\"uller flow on Teichm\"uller space (of, in general,  higher genus Riemann surfaces) to the genus one setting.   Here, the focus is on $\mathcal M$ being the Riemann moduli space of elliptic curves, or perhaps more to the point, its being the classifying space of flat lattices of (co-)area one.   Arnoux constructs a model of the unit tangent space of $\mathcal M$ by giving an explicit pair of lattice generators as well as a `flow direction' per lattice.    The flow direction corresponds to a  unit tangent vector whose basepoint on $\mathcal M$ represents  the class of the lattice.   The geodesic flow corresponds to squeezing a fundamental domain for the lattice in the vertical direction while stretching it in the horizontal direction.     Arnoux defines explicit coordinates and also finds the Gauss dynamical system as a factor of the dynamical system given by the first return by the geodesic flow to his explicit section.

 In  2013, we \cite{ArnouxSchmidtCross} used what the second author of this overview calls ``Arnoux's transversal" to show that each of the full one-parameter family of \cite{Nakada81}, known as the Nakada $\alpha$-continued fractions,  arises as the factor system of the first return by the geodesic flow to a corresponding cross section of the unit tangent bundle of $\mathcal M$.  The transversal arises from the identification of the unit tangent bundle of $\mathbb H$ with $\text{PSL}_2(\mathbb R)$, thus with the quotient of the real special linear group by its subgroup of nonzero constant multiples of the identity.   The identification allows for a map from $\mathbb R^2$ into this unit tangent bundle (images have  basepoints along horocycles).   There is a history of the use of planar dynamical systems to study continued fractions (again, see especially  \cite{Nakada81}); the above leads to a map from the domain of such a system into the unit tangent bundle, and in general this induces a measure preserving map from the system into the unit tangent bundle of $\mathcal M$.  However, this map usually fails to be injective. More important, the conjugate system may not agree with the {\em first} return map of the geodesic flow.  Still, for each $\alpha$-continued fraction  one finds that the geodesic flow   is a special flow over a shift which is essentially the natural extension of the $\alpha$-continued fraction transformation.   We mention that as noted in \cite{AbramsKatok}, the system of Adler and Flatto's linearized map is also the natural extension of an associated interval map.
 
 In 2017, Pohl \cite{Pohl} introduced a method to determine fundamental domains giving codings of geodesics that are particularly well suited to determining Maas forms by transfer operators (of the thermodynamic formalism).  For an extremely inviting introduction to these ideas, see her overview article with Zagier \cite{PohlZagier}.

  \section{ Subsequent developments}       The number of works exploring the connection between geodesic flows, cutting sequences and continued fractions continues to grow.       We hint at some particularly interesting directions by mentioning a limited selection of papers.    To help in making the selection, we   concentrated on works that cite \cite{Series}.

 \subsection{The modular surface setting}   One of the first to apply results of \cite{Series} was Pollicott \cite{Pollicott}, who in 1986  using the tools of Ruelle-Perron-Frobenius theory while taking the view of the geodesic flow on $\mathcal M$ as a special flow related to the continued fraction system   to reach  equidistribution results: The closed geodesics are equidistributed  by length with respect to the natural measure on the unit tangent bundle;  The quadratic irrationals are equidistributed in $[0,1]$ with respect to Gauss measure.   The latter result is deduced from the former by way of the relationship between closed geodesics and binary quadratic forms, a relationship clarified by Sarnak \cite{Sarnak}.
  
 In 2001, Grabiner and Lagarias \cite{GrabinerLagarias}  gave various results on cutting sequences associated to the standard fundamental domain of $\mathcal M$.    These included a response to a question of Adler-Flatto \cite{AdlerFlattoBull}: Cutting sequences can detect conformal properties.   In particular,   the symbolic dynamical shift system that \cite{GrabinerLagarias} define by   cutting sequences characterizes the standard fundamental domain of $\mathcal M$ up to isometries of $\mathbb H$ among all finite area hyperbolic polygons that are fundamental domains for a discrete subgroup of $\text{SL}_2(\mathbb R)$.

In a series of papers beginning with \cite{KatokUgarcoviciAbBegins} in 2010,   Katok and Ugarcovici (and later with Abrams) explored a two-dimensional family of continued fraction maps, all related to $\mathcal M$.   In \cite{KatokUgarcoviciApplAbAndFlo} they use cutting sequences, unit tangent vectors pointing into the standard fundamental domain of $\mathcal M$, and a reduction algorithm per member of their family to show   that the geodesic flow   is a special flow over a shift given in terms of their particular continued fraction map. 

 In 2012, Beardon, Hockman and  Short \cite{BeardonHockmanShort} used the Farey tesselation to introduce their notion of `geodesic continued fractions'.   In briefest terms, this is a method to describe general continued fractions related to $\mathcal M$ by way of following paths from $\infty$ to $x$ along the edges of the Farey tesselation.  

Since a large partial quotient of the continued fraction expansion indicates an excursion into the cusp of large order (in the sense of Humbert),  there is a close connection between  the `height' (depth  of penetration into the cusp at infinity) of a geodesic  on $\mathcal M$ and continued fraction expansions.   There are very many works using this to provide a geometric approach to questions in Diophantine approximation (already Ford's work is of this type).  In their 2017 \cite{BourgainKontorovich},  Bourgain and Kontorovich turn this around and, using cutting sequences, show the  existence of infinitely many low height closed `fundamental' geodesics. 
 
   Of the numerous papers studying properties of various continued fraction algorithms by way of cutting sequences, let us mention the 2018 \cite{BocaMerriman} work of Boca and Merriman  as well as    
     Merriman's  2023 \cite{Merriman}.

\subsection{The Farey map and additive continued fractions}

While it is standard to code numbers by their continued fraction expansion, using the infinite alphabet $\mathbb N$, the paper  \cite{Series} shows that there is a simpler coding, by the finite alphabet $\{L,R\}$, related to the Farey tree. One can indeed define a Farey map, topologically conjugate to the tent map, and study its properties.   The  Farey map  sends $x \in [0,1/2]$ to $x/(1-x)$ and $x \in [1/2, 1]$ to $(1-x)/x$.   There is an induction process that allows one to pass from the Farey map to the Gauss map, for this reason one can call the Farey map the {\em additive} version of the Gauss map.  

Given their close relation, one could expect that there is a coding of the geodesic flow on $\mathcal M$ in terms of the Farey map, and indeed  in 2019, Heersink \cite{Heersink} adjusted the cross section given by Series so as to obtain a cross section for the geodesic flow on the unit tangent bundle of $\mathcal M$ that projects to the Farey system.  

The properties of the Farey map are   quite different from those of the Gauss map:    the only absolutely continuous invariant measure has infinite mass, and the return time to the associated cross section has average zero, and thus the entropy of the Farey map is zero.    The Gauss map is related to the usual (multiplicative) Euclid algorithm on integers; in the same way, one easily proves that the Farey map is related to the additive Euclid algorithm, which is much slower.   This slowness is related to the fact that the natural invariant measure of the Farey map has infinite mass.

 A number of papers, see for example the 2021 paper of Cesarato and Vall\'ee \cite {CesarattoVallee}, have studied from a computer science viewpoint the differences between an additive and a multiplicative continued fraction, in various settings.

What are reasonably called the Farey convergents to a real $x$ form the set of the continued fraction convergents and the so-called {\em mediant} continued fraction convergents, values of the form $(a p_n + p_{n-1})/(a q_n + q_{n-1})$ with integers  $1\le a < a_n$.   Among the mediant and usual convergents are the best approximations of the first kind:  $p/q$ is such an approximation if $|x - p/q| < |x - c/d|$ for any rational $c/d$ with $d<q$.

 \subsection{Higher dimensional continued fractions}\label{ss:high}
 The regular continued fractions are in some sense Panglossian:   they converge well,   give one type of best approximation, and their additive version also gives a second version of best approximation.     For higher dimensional continued fractions algorithms, one must leave this best of all worlds,  for an overview of  some of this, see Berth\'e's \cite{Berthe} and the recent \cite{multiAuthorMultidim}.

 \subsection{Other surfaces}\label{ss:Other}    Series herself revisited the coding of geodesics in various manners.   Several of these are mentioned in a recent Notices of the AMS article on her work in hyperbolic geometry 
\cite{ParkerTanHypSeries}.   Among these are the co-edited collection \cite{BedfordKeaneSeries}, which provides an invitation to  the coding of geodesics and related matters.    Of Series' own research papers let us simply mention \cite{Series86}, where she applies cutting sequences and generalized continued fractions to  
obtain for each surface of constant negative curvature a representation of its geodesic flow as a special flow over  the shift involving certain words in a generating set of its fundamental group.  
Moreover, she clearly states and proves the relationship between the cutting sequence approach to the boundary expansions approach under weak restrictions (the fundamental domains considered are excluded from being triangles,  and must satisfy an {\em even corners} condition).      That is, \cite{Series86} is a deep and significant work building upon \cite{Series}.   

 In 1988, the aforementioned results on the Markoff spectrum of Diophantine approximation as related to simple closed geodesics on three covers of the modular surface were shown by Haas \cite{Haas} to hold for a large class of triples of hyperbolic surfaces.  The arguments rely in part on straightforward generalizations of \cite{Series}.
 
  The modular group, $\text{PSL}_2(\mathbb Z)$, is a member of the family of {\em Hecke (triangle Fuchsian) groups}.  This family is indexed by integers $q\ge 3$, with  each group given in terms of fractional linear transformations being generated  by $z \mapsto z +2 \cos \pi/q$ and $z\mapsto -1/z$.    The quotient of $\mathbb H$ by any of these groups is a genus zero surface with a cusp and two quotient singularities; $q=3$ gives $\mathcal M$.   The family is well-studied,  Rosen \cite{{Rosen}} gave continued fractions related to each group.    In 1986, Haas and Series \cite{HaasSeries} used cutting sequences based upon an analog of the Farey tesselation to answer questions about Diophantine approximation by elements of the group.     In 1996, Gr\"ochenig and Haas \cite{GroechnigHaas} considered `backwards' continued fractions in some generality, and coded geodesics in the unit tangent bundle of the surface related to the Farey tesselation for each $q$. 
  Planar natural extensions for the Rosen fractions were given by Nakada \cite{Nakada92} for the even case and by \cite{BurtonKraaikampSchmidt} in general.   As per \cite{ArnouxSchmidtCross}, these extensions determine cross-sections for the geodesic flow on the unit tangent bundle of the corresponding hyperbolic orbifold.    As part of  Mayer's long-term program of explicit applications of the Ruelle thermodynamic formalism (see his overview \cite{MayerBull}), with various co-authors he also gave codings of the geodesic flow on these surfaces; in  the 2008 \cite{MayerStroemberg} the coding explicitly uses unit tangent vectors directed inwards with respect to the analog of the standard fundamental domain.     
 
   In 1991, D.~Fried \cite{Fried} generalized the Artin approach to the setting of any finite index subgroup of a hyperbolic triangle  group with cusps.    The work focusses on representing conjugacy classes of elements (and thus coding of closed geodesics), gives the geodesic flow as a special flow, and follows Mayer's approach to applying the thermodynamic formalism.     Fried also makes clear why  dilogarithms so often arise in calculations involving volumes of unit tangent bundles of hyperbolic surfaces.

   Koebe \cite{Koebe} already considered various types of surfaces, including closed oriented surfaces of genus $g>1$.    In 1991, Adler and Flatto \cite{AdlerFlattoBull} applied their approach in this setting, insisting on a fundamental domain which is an $8g-4$-sided polygon (in an appendix they prove that such fundamental domains always exist) having the even corners property (which they call the `extension property').       Their work was revisited  in 2019 by Abrams and Katok \cite{AbramsKatok}, extending and simplifying the earlier work.  In particular, their notion of `short cycles' (a version of `matching') leads to new insights. 
   
   A compact closed Riemann surface has a unique hyperbolic structure, but many (let us say: nearly)  flat structures.   By integration, each nonzero holomorphic 1-form induces a {\em translation surface} structure on a Riemann surface.  The geodesic flow on the translation surface can be seen as a simple straight line flow along the surface.   Translation surfaces with a large  appropriately defined symmetry group are particularly susceptible to techniques related to hyperbolic surfaces.   For example, for the translation surface given by identifying opposite sides of a regular octagon by translations,   Smillie and Ulcigrai in their 2010 \cite{SmillieUlcigrai}  give cutting sequences in a manner reminiscent of \cite{Series}.   Questions arising from analogs with geodesic flow on hyperbolic surfaces continue to drive the theory of translation surfaces forward.     All of this is related to the notion of the Teichm\"uller flow, which we briefly mention below. 
 
 \subsection{Beyond real surfaces}   In 1994, Pollicott \cite{Pollicott2}  generalized the approach of \cite {Series}  to the quotient of hyperbolic 3-space by  $\text{SL}(2, \mathbb Z[i])$ so as to express  Selberg's zeta function for the geodesic flow here by way of the transfer operator of the appropriate boundary map.   Correspondingly, there is a complex continued fraction map.   The study of complex continued fractions goes back at least to the Hurwitz brothers;  Ford indicated how to use horoballs in place of his circles for this.   Results by way of the geodesic flow for  hyperbolic 3-manifolds are known, see for example the 2023 paper \cite{EiNakadaNatsui} of Ei, Nakada and Natsui.   However,  \cite{EiNakadaNatsui2}  shows that there are quite reasonable complex continued fractions for which this method is not fully appropriate. 
 
 Lukyanenko and Vandehey have generalized \cite{Series} in various ways, see for example their 2023 \cite{LukyanenkoVandehey}, where they define certain continued fraction maps in various dimensions for which they can use a geodesic flow to study the map.  (In that paper, they in fact take the boundary expansion approach, but to obtain good codings there is a jump past certain `small digits'.)

  The notions of geodesic flow can be made sense of in yet other `non-real' settings.     In his Ph.D. dissertation Artin  introduced polynomial continued fractions for Laurent series over finite fields.  In 2007 Broise-Alamichel and Paulin \cite{BroisePaulin} in a sense  mimicked \cite{Series} to show that  the first return map on a certain 
cross section to a geodesic flow on the appropriate quotient of the Bruhat-Tits
tree for  $\text{PGL}_2$ of the Laurent series field gives Artin's polynomial map. 
 
 \subsection{Other flows}\label{ss:OtherFlows} 
Besides following the flow of unit tangent vectors along a geodesic, one can follow the flow along a horocycle.   The study of the horocycle flow on surfaces has a long history.    A  coding for the horocycle flow on the unit tangent bundle of $\mathcal M$ using an interval map related to the Farey fractions was given in 2013 by Athreya and Cheung \cite{AthreyaCheung}.    

The aforementioned Arnoux coding approach is based upon the well known identification of the geodesic flow on the unit tangent bundle of $\mathcal M$ with a natural flow on area one lattices in $\mathbb R^2$.   This flow directly generalizes to higher dimensions:  the space of area one lattices in $\mathbb R^n$ can be identified with $\text{SL}_n(\mathbb Z)\backslash \text{SL}_n(\mathbb R)$ and there is an action of diagonal real matrices giving the so-called diagonal flow on  this space.    There are many important works on this flow (as well as on a great variety of generalizations) and especially on its connection to Diophantine approximation, see for example the 1998 work of Kleinbock and Margulis \cite{KleinbockMargulis}.    Returning to the classical setting, in 2007  Ghys  (again see \cite{Ghys} for a discussion) using classical means to identify  the unit tangent bundle of $\mathcal M$ with the complement of the trefoil knot in the sphere $S^3$ and the diagonal flow identified the set of knots formed by the set of  closed geodesics of $\mathcal M$.

The closed geodesics of $\mathcal M$ correspond to $\text{SL}_2(\mathbb Z)$-conjugacy classes of matrices with trace of absolute value greater than 2.  Already Poincar\'e used the fact that such a matrix defines a linear automorphism of the square torus,  with expanding and contracting directions corresponding to the distinct eigenspaces of the matrix.   Directly generalizing this are Anosov diffeomorphisms and their flows.     
Also related is the notion of a {\em pseudo-Anosov} homeomorphism of a real surface; these were defined by W.~Thurston, who showed how to view such a homeomorphism as being given appropriately by a real $2\times 2$ matrix of distinct eigenspaces.  The mapping cylinder of a pseudo-Anosov has a suspension flow on it with many interesting properties, see say \cite{McMullenFibered}.   

The Teichm\"uller flow briefly mentioned above is also an area of intense study,  for its connections to the aforementioned  translation surfaces setting, see say \cite{McMullenBill} and \cite{Zorich}. Let us just say that the geodesic flow on the modular surface can be seen as a  {\em ``toy model"}  for the Teichm\"uller flow, and that the starting point of the recent decades of progress in the study of translation surfaces and the Teichm\"uller flow came from associating a continued fraction (the Rauzy induction of interval exchange transformations) to the action of the Teichmüller flow: the essential result of the 1982 paper by Veech \cite{Veech} consists in building an explicit natural extension for Rauzy induction, an invariant measure for a suspension, and proving that this measure is finite.

Finally, let us remark that one can start from a given continued fraction (such as the Brun or Selmer continued fraction) and try, as was done for Rauzy induction, to build a natural extension, an invariant measure, and a flow on a modular manifold, associated with a form of induction on a family of dynamical systems.  For example, this program has been  partially realized in the recent paper \cite{multiAuthorMultidim} for algorithms satisfying a generalized Pisot property; however, the modular space and the flow built in this way have not yet been identified with other classical objects, as is the case for the modular surface.

 \section{Conclusion:    Highly recommended reading}  It is hard to find better words than those of P.~J.~Nicholls in his review of Series' \cite{Series}, where he wrote 
 {\it  This paper is beautifully written and the results are very elegant. Both this paper and the author's recent article [Math. Intelligencer 7 (1985), no. 3, 20–29; MR0795536] should be required reading for those interested in the symbolic dynamics of geodesic flows. }     
 
 Besides that cited by Nicholls, we heartily recommend Series' earlier Intelligencer article \cite{SeriesIntell82}.   Let us again point to \cite{BedfordKeaneSeries}  as a fine  introduction to the general area of \cite{Series}; for a book aimed at a wider audience see \cite{IndraPs}.   Another fine introductory textbook to geodesic (and horocyclic) flows is Dal'bo's \cite{Dal'Bo}.  Excellent textbook introductions to continued fractions abound --- with the classic treatment being that  of Hardy and Wright \cite{HardyWright}.   Another classic reference is Perron's two volume \cite{Perron1, Perron2}, which gives state-of-the-art, as of the mid-20th century, coverage of the theory of continued fractions.  Continued fractions as an example of ergodic theory is elegantly given in Billingsley's textbook \cite{Billingsley}.   Farey fractions play a recurring role in Hatcher's  \cite{Hatcher}   topological take on number theory.


Finally, both the list above and the citations throughout are necessarily selective, and many excellent contributions related to \cite{Series} had to be omitted.


\begin{thebibliography}{ABMST}  


 
\bibitem[AK]{AbramsKatok} A. Abrams and S. Katok, {\em
Adler and Flatto revisited: cross-sections for geodesic flow on compact surfaces of constant negative curvature}, 
Studia Math. 246 (2019), no. 2, 167--202. 

\bibitem[AF]{AdlerFlatto} R. Adler and L. Flatto, {\em Cross section maps for geodesic flows. I. The modular surface}, in  Ergodic theory and dynamical systems, II (College Park, Md., 1979/1980), pp. 103--161, Progr. Math., 21, Birkhäuser, Boston, MA, 1982.


\bibitem[AF2]{AdlerFlattoBack}
\bysame, {\em The backward continued fraction map and geodesic flow},  Ergodic Theory Dynam. Systems 4 (1984), no. 4, 487--492.


\bibitem[AF3]{AdlerFlattoBull}
\bysame, {\em Geodesic flows, interval maps and symbolic dynamics},  Bull. Amer. Math. Soc., 25 (1991), No. 2, 229-334.


\bibitem[Ar]{ArnouxCodage} P.~Arnoux, {\em Le codage du flot g\'eod\'esique sur la surface modulaire. } Enseign. Math. (2) 40 (1994), no. 1-2, 29--48.

\bibitem[ABMST]{multiAuthorMultidim} P.~ Arnoux,V.~Berth\'e, M.~Minervino, W.~Steiner, J.~M.~Thuswaldner, {\em 
Nonstationary Markov partitions and multidimensional continued fraction algorithms},   
 \href{https://arxiv.org/abs/2508.16441}{preprint 2025, 147 pp.: Arxiv/2508.16441}. 

\bibitem[AS]{ArnouxSchmidtCross} P. Arnoux and T. A. Schmidt, {\em Cross sections for geodesic flows and $\alpha$-continued fractions},  Nonlinearity 26 (2013), 711--726.


	
\bibitem[Art]{Artin} E. Artin, {\em Ein mechanisches System mit quasi-ergodischen Bahnen}, Abh. Math. Sem. Hamburg 3
(1924)
170--175 (and {\em Collected Papers}, Springer-Verlag, New
York,
1982,  499--505).  


\bibitem[AC]{AthreyaCheung} J.~S.~Athreya and Y.~Cheung,  {\em A Poincar\'e section for the horocycle flow on the space of lattices}, Int. Math. Res. Not. IMRN 2014, no. 10, 2643--2690.

\bibitem[BHS]{BeardonHockmanShort} A.~F.~Beardon,  M.~Hockman,  and I.~Short,  
{\em Geodesic continued fractions},
Michigan Math. J. 61 (2012), no. 1, 133--150. 

\bibitem[BKS]{BedfordKeaneSeries}   Ergodic Theory, symbolic dynamics, and hyperbolic spaces,  T.~Bedford, M.~Keane, C.~Series, eds., Oxford Univ. Press, 1991. 

\bibitem[B]{Berthe}  V.~Berth\'e,  
{\em Multidimensional Euclidean algorithms, numeration and substitutions}, 
Integers 11B (2011), Paper No. A2, 34 pp. 

  
 
\bibitem[Bi]{Billingsley} P.~ Billingsley, {\em Ergodic Theory and Information}, Wiley, 1965.

\bibitem[BoM]{BocaMerriman} F.~Boca and C.~Merriman,  
{\em Coding of geodesics on some modular surfaces and applications to odd and even continued fractions}, Indag. Math. (N.S.) 29 (2018), no. 5, 1214--1234. 

\bibitem[BK]{BourgainKontorovich} J.~Bourgain and A.~ Kontorovich,  
{\em Beyond expansion II: low-lying fundamental geodesics.},
J. Eur. Math. Soc. (JEMS) 19 (2017), no. 5, 1331--1359. 

\bibitem[BP]{BroisePaulin}    A.~Broise-Alamichel and F.~Paulin, {\em Dynamiques sur le rayon modulaire et fractions continues en caractéristique $p$},  [Dynamics on the modular ray and continued fractions in characteristic p] J. Lond. Math. Soc. (2) 76 (2007), no. 2, 399--418.

\bibitem[BuKS]{BurtonKraaikampSchmidt} R.~M.~Burton, C.~Kraaikamp and T.~A.~Schmidt, {\em Natural extensions for the Rosen fractions},
Trans. Amer. Math. Soc., 352 (2000), 1277--1298.

\bibitem[C]{Conway}  J.~H.~Conway, 
{\em An enumeration of knots and links, and some of their algebraic properties.} In: Computational Problems in Abstract Algebra (Proc. Conf., Oxford, 1967), pp. 329--358, Pergamon, Oxford-New York-Toronto, Ont., 1970


\bibitem[CN]{CovenNitecki}   E.~M.~Coven and Z.~Nitecki, 
{\em On the genesis of symbolic dynamics as we know it},
Colloq. Math. 110 (2008), no. 2, 227--242.

\bibitem[CV]{CesarattoVallee} E. Cesaratto and B. Vall\'{e}e,
 {\em Gaussian behavior of quadratic irrationals},
Acta Arith, 197 {2021},
   no. 2, 159--205,
	

\bibitem[D]{Dedekind}    R. Dedekind, {\em Schreiben an Herrn Borchardt \"uber die Theorie der elliptischen Modulfunktionen}, J. Reine Angew. Math. 83 (1877) 265--292.

\bibitem[Da]{Dal'Bo}  F.~Dal'Bo,   Geodesic and horocyclic trajectories.  Universitext. Springer-Verlag London, 2011.

\bibitem[ENN]{EiNakadaNatsui}  H.~Ei,  H.~Nakada and R.~Natsui, 
{\em On the ergodic theory of maps associated with the nearest integer complex continued fractions over imaginary quadratic fields}, 
Discrete Contin. Dyn. Syst. 43 (2023), no. 11, 3883--3924. 

\bibitem[ENN2]{EiNakadaNatsui2} \bysame, 
{\em On the dynamics of a complex continued fraction map which contains the Gauss map as its real number section}, 
Adv. Math. 472 (2025), Paper No. 110286, 36 pp. 

\bibitem[F]{Ford1938} L.~R.~Ford, {\em Fractions},  The American Mathematical Monthly, Vol. 45, No. 9 (Nov., 1938), pp. 586--601.

\bibitem[Fr]{Fried} D.~Fried, {\em Symbolic dynamics for triangle groups}, Invent. Math., 125 (1996), 487--521.

\bibitem[G]{Ghys} E.~Ghys, {\em G\'eod\'esiques sur les surfaces \`a courbure n\'egative}, Le\c{c}ons de math\'ematiques d’aujourd’hui, volume 4, pr\'esent\'e par Fr\'ed\'eric Bayart et Eric Charpentier, Belin (2009)
 
\bibitem[GL]{GrabinerLagarias} D.~J.~Grabiner and J.~C.~Lagarias,  
{\em Cutting sequences for geodesic flow on the modular surface and continued fractions}, 
Monatsh. Math. 133 (2001), no. 4, 295–339. 

\bibitem[GH]{GroechnigHaas} K.~Gr\"ochenig and A.~Haas, {\em Backward continued fractions, Hecke groups and invariant measures
for transformations of the interval}, Ergodic Theory Dynam. Systems, 16 (1996), 1241--1274.
 
\bibitem[Ha]{Haas} A.~Haas,
{\em  Diophantine approximation on hyperbolic orbifolds},
Duke Math. J. 56 (1988), no. 3, 531--547. 

\bibitem[HS]{HaasSeries} A.~Haas and C.~Series, 
{\em The Hurwitz constant and Diophantine approximation on Hecke groups},
J. London Math. Soc. (2) 34 (1986), no. 2, 219--234. 

\bibitem[H]{Hadamard} J.~Hadamard, {\em Les surfaces \`a courbures oppos\'ees et leurs lignes g\'eod\'esiques}, J. Math. Pures
Appl. 4 (1898), p. 27--74. OEuvres, tome II, p. 729--775.

 
\bibitem[HW]{HardyWright} 
G. H.~Hardy and E.~M.~Wright, An introduction to the theory of numbers. Sixth edition. Revised by D. R. Heath-Brown and J. H. Silverman. With a foreword by Andrew Wiles. Oxford University Press, Oxford, 2008.

\bibitem[Ha]{Hatcher} 
A.~Hatcher,  Topology of numbers. American Mathematical Society, Providence, RI, 2022. 
 

\bibitem[He]{HedlundMod}  G.~A.~Hedlund, {\em A metrically transitive group defined by the modular group},  Amer. J. Math. 57
(1935), 668--678.

\bibitem[He2]{HedlundClosedSurfs}
\bysame,  {\em On the metrical transitivity of geodesies on closed surfaces of constant negative
curvature},  Ann. Math. 35 (1934), 787--808.

\bibitem[Hi]{HirzebruchQuotient}
F.~Hirzebruch, {\em \"Uber vierdimensionale Riemannsche Fl\"achen mehrdeutiger analytischer
Funktionen von zwei komplexen Ver\"anderlichen}, Math. Ann. 126 (1953), 1--22.
 
\bibitem[Hi2]{HirzebruchCusp} 
 {\em The Hilbert modular group, resolution of the singularities at the cusps and related problems}, in: Séminaire Bourbaki, 23\`eme ann\'ee (1970/1971), Exp. No. 396, pp. 275--288, Lecture Notes in Math., Vol. 244, Springer, Berlin-New York, 1971. 

 
\bibitem[Hee]{Heersink} B.~Heersink, 
{\em Distribution of the periodic points of the Farey map},
Comm. Math. Phys. 365 (2019), no. 3, 971--1003. 

\bibitem[Ho]{Hopf} E.~Hopf,
{\em Ergodic theory and the geodesic flow on surfaces of constant negative curvature},  
Bull. Amer. Math. Soc. 77,  863--877 (1971).

 
\bibitem[Hu]{Humbert}  G.~Humbert, 
{\em Sur les fractions continues ordinaires et les formes
quadratiques binaires ind\'efinies}, 
Journal de mathématiques pures et appliquées 7e s\'erie, tome 2 (1916), p. 104--154.

\bibitem[K]{Koebe}  P.~Koebe, {\em Riemannsche Mannigfaltigkeiten und nicht euklidische Raumformen}, Sitzungsberichte
der Preu{\ss}ischen Akademie der Wissenschaften, I (1927),414--457.
 
\bibitem[KU]{KUbAMS} 
S.~Katok, I.~Ugarcovici, 
{\em Symbolic dynamics for the modular surface and beyond},
Bull. Amer. Math. Soc. (N.S.) 44 (2007), no. 1, 87--132.  
 
\bibitem[KU2]{KatokUgarcoviciAbBegins} \bysame, {\em Theory of (a,b)-continued fraction transformations and applications}, Electron. Res. Announc. Math. Sci. 17 (2010), 20–33. 

\bibitem[KU3]{KatokUgarcoviciApplAbAndFlo} \bysame, {\em Applications of $(a,b)$-continued fraction transformations},
Ergodic Theory Dynam. Systems 32 (2012), no. 2, 755--777. 

\bibitem[KM]{KleinbockMargulis}   D.~Y.~Kleinbock and  G.~A.~Margulis, 
{\em Flows on homogeneous spaces and Diophantine approximation on manifolds}, 
Ann. of Math. (2) 148 (1998), no. 1, 339--360.

\bibitem[LV]{LukyanenkoVandehey}  A.~Lukyanenko and J.~Vandehey,  
{\em Ergodicity of Iwasawa continued fractions via markable hyperbolic geodesics},
Ergodic Theory Dynam. Systems 43 (2023), no. 5, 1666--1711. 
 
\bibitem[M]{MayerBull}  D.~Mayer, {\em The thermodynamic formalism approach to Selberg’s zeta function for
$\text{PSL}(2, Z)$}, Bull. Amer. Math. Soc. (N.S.) 25 (1991), 55--60.


\bibitem[MS]{MayerStroemberg}  D.~Mayer and F.~Str\"omberg, {\em Symbolic dynamics for the Geodesic flow on Hecke
surfaces}, Journal of Modern Dynamics 2 (2008), 581--627.


\bibitem[Mc]{McMullenFibered}  C.~T.~McMullen, {\em Polynomial invariants for fibered 3-manifolds and Teichm\"uller geodesics for foliations},  Ann. Sci. \'Ecole Norm. Sup. (4) 33 (2000), no. 4, 519--560.

\bibitem[Mc2]{McMullenBill}  \bysame, {\em Billiards and Teichm\"uller curves},  Bull. Amer. Math. Soc. (N.S.) 60 (2023), no. 2, 195--250. 

\bibitem[Me]{Merriman} C.~ Merriman, 
{\em Geodesic flows and the mother of all continued fractions},
Int. J. Number Theory 18 (2022), no. 4, 931--953.

\bibitem[Mi]{Minkowski} H.~ Minkowski, 
{\em G\'en\'eralisation de la th\'eorie des fractions continues},
Annales scientifiques de l’\'E.N.S. 3e série, tome 13 (1896), pp. 41--60. 

\bibitem[Moe]{Moeckel} R.~Moeckel,  
{\em Geodesics on modular surfaces and continued fractions}, 
Ergodic Theory Dynam. Systems 2 (1982), no. 1, 69--83.  

\bibitem[Mor]{Morse1}  M.~Morse, {\em  A one-to-one representation of geodesics on a surface of negative curvature},
Amer J Math 43  (1921) 33--51.

\bibitem[Mor2]{Morse2}  \bysame, {\em  Recurrent geodesics on a surface of negative curvature},  Trans. Amer. Math. Soc. XXII (1921), 84--100.

\bibitem[Mor3]{MorseNotes}  \bysame, {\em Symbolic dynamics},  Institute for Advanced Study Notes, Princeton (1966) (unpublished). (First written 1938.)
 
\bibitem[MSW]{IndraPs} 
  D.~Mumford, C.~Series, and D.~Wright, Indra's pearls. The vision of Felix Klein. Cambridge University Press, New York, 2002.
 
  
\bibitem[N81]{Nakada81}  H.~Nakada,
{\em Metrical theory for a class of continued fraction transformations and their natural extensions},
Tokyo J. Math., 4,   399--426 (1981).


\bibitem[N92]{Nakada92} \bysame,
{\em Continued fractions, geodesic flows and Ford circles},
In:   Algorithms, fractals, and dynamics: Okayama/Kyoto, 1992,
 (1995) 179--191.
 
\bibitem[Ni]{Nielsen} J.~Nielsen, {\em Untersuchungen zur Topologie der geschlossenen zweiseitigen Fl\"achen},  Acta
Math 50, (1927), 189--358.

\bibitem[PT]{ParkerTanHypSeries}  J. R.~Parker and S.~P.~Tan, {\em Caroline Series and hyperbolic geometry}, Notices Amer. Math. Soc. 70 (2023), no. 3, 380--389.  

\bibitem[P1]{Perron1} O.~Perron,   Die Lehre von den Kettenbr\"uchen. Bd I. Elementare Kettenbrüche. (German) 3te Aufl. B. G. Teubner Verlagsgesellschaft, Stuttgart, 1954. 
 
 
\bibitem[P2]{Perron2} \bysame,  Die Lehre von den Kettenbr\"uchen. Dritte, verbesserte und erweiterte Aufl. Bd. II. Analytisch-funktionentheoretische Kettenbrüche. (German) B. G. Teubner Verlagsgesellschaft, Stuttgart, 1957.  

\bibitem[Po]{Pohl}  A.~Pohl, {\em Symbolic dynamics for the geodesic flow on two-dimensional hyperbolic good orbifolds},  Discrete Contin. Dyn. Syst. 34 (2014), no. 5, 2173--2241.

\bibitem[PZ]{PohlZagier} A.~Pohl and D.~Zagier, 
{\em Dynamics of geodesics, and Maass cusp forms},
Enseign. Math. 66 (2020), no. 3-4, 305–340. 
 
\bibitem[Pol]{Pollicott} M.~Pollicott,  
{\em Distribution of closed geodesics on the modular surface and quadratic irrationals},
Bull. Soc. Math. France 114 (1986), no. 4, 431--446. 


\bibitem[Pol]{Pollicott2} \bysame, 
{\em The Picard group, closed geodesics and zeta functions},
Trans. Amer. Math. Soc. 344 (1994), no. 2, 857--872.

\bibitem[R]{Rosen}   D.~Rosen, {\em A class of continued fractions associated with certain properly
discontinuous groups}, Duke Mathematical Journal 21 (1954), 549--563.

\bibitem[S]{Sarnak}  P.~Sarnak, 
{\em Class numbers of indefinite binary quadratic forms}, J. Number Theory 
Vol. 15, 1982, pp. 229--247.

\bibitem[S81]{Series81} C.~Series, 
{\em  On coding geodesics with continued fractions}, Ergodic theory (Sem., Les Plans-sur-Bex, 1980) (French), pp. 67--76, Monogr. Enseign. Math., 29, Univ. Genève, Geneva, 1981.


\bibitem[S]{Series} C. Series, {\em The modular surface and continued fractions}, J. Lond. Math. Soc. (2), {\bf 31}, no.~1, (1985), 69-80.

\bibitem[S82]{SeriesIntell82}
\bysame, {\em Non-Euclidean geometry, continued fractions, and ergodic theory}, Math. Intelligencer 4 (1982), no. 1, 24--31.
 
\bibitem[S85]{SeriesGeomMarkoffNumbers}
\bysame, {\em The geometry of Markoff numbers},
Math. Intelligencer 7 (1985), no. 3, 20--29.

\bibitem[S86]{Series86}
\bysame, {\em Geometrical Markov coding of geodesics on surfaces of constant negative curvature}, Ergodic Theory Dynam. Systems 6 (1986), no. 4, 601--625.





\bibitem[S91]{SeriesArtInBook}
\bysame, {\em Geometrical methods of symbolic coding}.  In \cite{BedfordKeaneSeries}, pp. 125--151.


 
\bibitem[SU]{SmillieUlcigrai} J.~Smillie and C.~Ulcigrai,  
{\em Geodesic flow on the Teichmüller disk of the regular octagon, cutting sequences and octagon continued fractions maps}, In:  Dynamical numbers—interplay between dynamical systems and number theory, 29--65,
Contemp. Math., 532, Amer. Math. Soc., Providence, RI, 2010. 

\bibitem[T]{Thurston}  W.~P.~Thurston,  
{\em On proof and progress in mathematics},
Bull. Amer. Math. Soc. (N.S.) 30 (1994), no. 2, 161--177.



\bibitem[Ve]{Veech} W.A. Veech,
{\em Gauss measures for transformations on the space of interval exchange maps},
Ann. of Math. 115 (1982) 201--242

\bibitem[Vu]{Vulakh} L.~Ya.~Vulakh,  {\em The Markov spectra for Fuchsian groups}, 
Trans. Amer. Math. Soc. 352 (2000), no. 9, 4067--4094. 

\bibitem[Z]{Zorich}
A.~ Zorich, 
{\em Flat surfaces}, in Frontiers in number theory, physics, and geometry. I, 437–583, Springer, Berlin, 2006. 

\end{thebibliography}
\end{document}